\documentclass[12pt,english,frenchb]{article}

\usepackage[T1]{fontenc}
\usepackage[latin1]{inputenc}
\usepackage{amsthm,amssymb,amsmath,mathrsfs}
\usepackage{babel}

\newcommand{\noun}[1]{\textsc{#1}}

 \NoAutoSpaceBeforeFDP

 \theoremstyle{plain}    
 \newtheorem{thm}{Théorème}[section]
 \theoremstyle{plain}    
 \newtheorem*{thm*}{Théorème} 
 \theoremstyle{plain}    
 \newtheorem{cor}[thm]{Corollaire} 
 \theoremstyle{plain}    
 \newtheorem{lem}[thm]{Lemme} 
 \theoremstyle{definition}
  \newtheorem{example}[thm]{Exemple}
 \theoremstyle{remark}
 \newtheorem{rem}[thm]{Remarque}

\numberwithin{equation}{section}

\renewcommand{\mathcal}{\mathscr}
\renewcommand{\mathbf}{\mathbb}
\newcommand{\im}{\mathrm{Im}\,}

\renewcommand{\leq}{\leqslant}
\renewcommand{\geq}{\geqslant}

\begin{document}

\title{Sur des variétés de Cauchy-Riemann dont la forme de Levi a une valeur
propre positive}

\date{}
\author{Olivier Biquard%
\thanks{The author is a member of EDGE, Research Training Network HPRN-CT-2000-00101,
supported by the European Human Potential Programme.
}}

\maketitle
\begin{abstract}
\selectlanguage{english}
For certain real hypersurfaces in the projective space, of signature
\( (1,n) \), we study the filling problem for small deformations
of the CR structure (the other signatures being well understood).
We characterize the deformations which are fillable, and prove that
they have infinite codimension in the set of all CR structures. This
result generalizes the cases of the 3-sphere and of signature \( (1,1) \)
to higher dimension.
\selectlanguage{frenchb}
\end{abstract}
Un problème classique d'analyse complexe consiste à savoir quand une
variété CR est le bord d'une variété complexe. Concentrons nous sur
la question, déjà intéressante, des petites déformations. Plus précisément,
soit \( N \) une variété complexe, telle que la forme de Levi du
bord \( \partial N \) soit non dégénérée, nous essayons de répondre
aux deux problèmes suivants :

\begin{enumerate}
\item \emph{problème de remplissage :} une déformation CR de \( \partial N \)
est-elle le bord d'une déformation complexe de \( N \) ? dans l'affirmative,
on dira que la déformation de \( \partial N \) est \emph{remplissable}
;
\item \emph{problème de stabilité :} supposons que \( N \) soit un ouvert
dans une variété complexe compacte \( V \), est-ce qu'une déformation
remplissable de \( \partial N \) est induite par une déformation
de l'hypersurface \( \partial N \) dans \( V \) ?
\end{enumerate}
Répétons que dans cet article, nous ne regardons que les petites déformations,
donc, par exemple dans le problème de remplissage, nous ne cherchons
que les remplissages proches de la variété initiale \( N \) ; cela
laisse de côté le problème global de remplissage de \( \partial N \)
par n'importe quelle variété complexe à bord.

Deux belles applications du théorème de Nash-Moser permettent de résoudre
les deux problèmes dans beaucoup de cas. En effet le premier problème
est résolu par le théorème suivant \cite{Kir79}.

\begin{thm*}
[Kiremidjian]Si la forme de Levi de \( \partial N \) a 0 ou au moins
2 valeurs propres strictement positives, et si \( H^{2}_{c}(N,T^{1,0}_{N})=0 \),
alors toutes les petites déformations de \( \partial N \) sont remplissables.
\end{thm*}
Le second problème est résolu par le théorème suivant \cite{Ham77}.

\begin{thm*}
[Hamilton]Si la forme de Levi de \( \partial N \) a 0 ou au moins
2 valeurs propres strictement négatives, et si \( H^{1}(N,T^{1,0}_{N})=0 \),
alors les déformations complexes de \( N \) sont induites par une
déformation de \( \partial N \) dans \( V \).
\end{thm*}
Supposons que dans une variété complexe compacte \( V \) on ait une
hypersurface réelle \( X \), séparant \( V \) en deux morceaux,
\( V=N\cup M \) et \( \partial N=\partial M=X \), alors les formes
de Levi de \( N \) et \( M \) ont des signatures opposées, donc
la condition sur la forme de Levi de \( N \) dans le théorème de
Kiremidjian est équivalente à la condition sur la forme de Levi de
\( M \) dans le théorème de Hamilton. Cela suggère que le problème
de remplissage pour \( N \) pourrait être lié au problème de stabilité
pour \( M \), ce qui est le cas, comme on le verra plus loin.

Les deux théorèmes ci-dessus laissent ouvert le cas où la forme de
Levi, non dégénérée, a exactement une valeur propre positive ou négative.
Dans ce papier, nous allons nous concentrer sur l'exemple le plus
simple de cette situation.

\begin{example}
\label{ex-base}L'hypersurface réelle \begin{equation}
\label{def-X}
X^{2n-1}=\{|z_{1}|^{2}+|z_{2}|^{2}=|z_{3}|^{2}+\cdots +|z_{n+1}|^{2}\}\subset P^{n}
\end{equation}
est le bord des deux domaines complexes de \( P^{n} \) définis par
\begin{eqnarray}
N^{n} & = & \{|z_{1}|^{2}+|z_{2}|^{2}<|z_{3}|^{2}+\cdots +|z_{n+1}|^{2}\},\label{def-N} \\
M^{n} & = & \{|z_{1}|^{2}+|z_{2}|^{2}>|z_{3}|^{2}+\cdots +|z_{n+1}|^{2}\},\label{def-M} 
\end{eqnarray}
 dont les formes de Levi ont pour signatures respectives \( (1,n-2) \)
et \( (n-2,1) \).
\end{example}
Les théorèmes de Kiremidjian et Hamilton nous indiquent déjà qu'on
a toujours une réponse positive aux problèmes de remplissage de \( X^{2n-1} \)
du côté de \( M^{n} \) pour \( n\neq 3 \) d'une part, de stabilité
de \( X^{2n-1} \) du côté de \( N^{n} \) pour \( n>3 \) d'autre
part.

Pour \( n=2 \), on a \( X=S^{3} \), avec \( N=B^{4} \), et la réponse
au problème de déformation est négative (il y a un espace de dimension
infinie d'obstruction), tandis que la réponse au problème de stabilité
de \( N \) est positive, voir les travaux \cite{BurEps90,Eps92,Lem92,Bla94}.

Pour \( n=3 \), donc en signature (1,1), j'ai montré \cite{BiqADE}
que les réponses aux problèmes de déformation et de stabilité sont
toujours négatives (et à chaque fois, l'obstruction est de dimension
infinie).

Le but de cette note est d'étendre ces résultats au cas des \( X^{2n-1} \)
en dimension supérieure (\( n>3 \)). Avant d'énoncer les résultats,
observons qu'on a une action de cercle sur \( X^{2n-1} \) donnée
par (\( \zeta \in S^{1} \)) \[
[z_{1}:\cdots :z_{n+1}]\longrightarrow [\zeta z_{1}:\zeta z_{2}:z_{3}:\cdots :z_{n+1}],\]
 qui s'étend en une action du disque \( \Delta  \) sur \( N^{n} \).
Cette action de \( S^{1} \) permet de décomposer les tenseurs sur
\( X \) en séries de Fourier, et en particulier les déformations
CR de \( X \), qui sont paramétrées par des (0,1)-formes sur \( X \)
à valeurs dans \( T^{1,0}_{X} \).

\begin{thm}
\label{th-def}Pour tous les \( n\geq 2 \), la réponse au problème
de remplissage de \( X^{2n-1} \) comme bord de \( N^{n} \) est négative,
et l'espace d'obstruction est de dimension infinie. En outre, une
petite déformation CR de \( X \) est remplissable si et seulement
si ses coefficients de Fourier strictement négatifs s'annulent (modulo
les contactomorphismes de \( X \)).
\end{thm}
Plus précisément, on peut trouver une jauge par rapport à l'action
des contactomorphismes de \( X \) telle que la structure CR soit
remplissable si et seulement si ses coefficients de Fourier strictement
négatifs s'annulent dans cette jauge. L'espace des obstructions est
explicitement calculé dans le lemme \ref{lem-obstructions}. En dimensions
3 et 5, le théorème est bien sûr déjà connu, comme mentionné plus
haut.

\begin{cor}
\label{th-stab}Pour tous les \( n\geq 2 \), la réponse au problème
de stabilité de \( M^{n}\subset P^{n} \) est négative, et l'espace
d'obstruction est de dimension infinie. En outre, pour \( n\neq 3 \),
une petite déformation CR de \( X^{2n-1} \) est remplissable du côté
de \( N^{n} \) si et seulement si son remplissage du côté de \( M^{n} \)
est stable.
\end{cor}
Là encore, la situation en dimensions 3 et 5 était déjà connue. À
noter qu'en dimension 5, il est de plus montré dans \cite{BiqADE}
que les structure CR obtenues par déformation de \( X \) dans \( P^{3} \)
sont de codimension infinie dans l'espace des structures CR remplissables
d'un côté.

Enfin, nous nous contentons ici d'étudier la série d'exemples (\ref{ex-base}),
mais il semble plausible que les phénomènes observés dans le théorème
\ref{th-def} et le corollaire \ref{th-stab} restent valables en
général sur les variétés CR de signature \( (1,n-2) \), avec éventuellement
des obstructions supplémentaires de dimension finie.

La méthode utilisée s'appuie fortement sur mon article \cite{BiqADE}
; en effet, bien qu'il soit consacré à la dimension 5, certains résultats
y sont démontrés en dimension quelconque. En particulier, la paramétrisation
du groupe des contactomorphismes y est réalisée, et démontrée la nécessité
que les coefficients de Fourier soient positifs pour avoir un remplissage.
Pour obtenir le cas de la dimension supérieure, il faut obtenir une
bonne paramétrisation des structures CR (section \ref{sec-1}), et
calculer l'action du groupe des contactomorphismes (section \ref{sec-2})
; cela est réalisé ici par des techniques différentes de celles de
\cite{BiqADE}, s'étendant plus facilement en dimension plus grande.

\section{\label{sec-1}La variété des structures CR sur \protect\( X\protect \)}

Considérons l'hypersurface \( X \) de \( P^{n} \) définie dans l'exemple
\ref{ex-base}. Remarquons que la projection \( [z_{1}:\cdots :z_{k+1}]\rightarrow ([z_{1}:z_{2}],[z_{3}:\cdots :z_{n+1}]) \)
permet de voir \( X \) comme l'espace total du fibré en cercles \( \mathcal{O}(-1,1) \)
au-dessus de \( P^{1}\times P^{n-2} \).

\subsubsection*{Intégrabilité}

Nous noterons \( H\subset TX \) la distribution de contact induite
par la structure CR. Puisqu'une déformation de la structure de contact
reste équivalente à \( H \), on peut supposer qu'une déformation
CR de \( X \) préserve \( H \), et est donc paramétrée par un tenseur
\( \phi \in \Omega ^{0,1}_{X}\otimes T^{1,0}_{X} \), de sorte que
l'espace \( T^{0,1} \) de la nouvelle structure complexe soit constitué
des \[
h+\phi _{h},\quad h\in T^{0,1}_{X}.\]
 Dans ce cas, la nouvelle structure complexe est intégrable à condition
que \begin{equation}
\label{con-int}
\overline{\partial }\phi +\frac{1}{2}[\phi ,\phi ]=0,
\end{equation}
 où \( \overline{\partial } \) est l'opérateur naturel du fibré holomorphe
\( T'X=T^{\mathbf{C}}X/T^{0,1}X \) sur \( X \), et le crochet \( [\phi ,\phi ]\in \Omega ^{0,2}_{X}\otimes T^{1,0}_{X} \)
est défini par \[
[\phi ,\psi ]_{h,k}=[\phi _{h},\psi _{k}]-\phi _{[h,\psi _{k}]^{0,1}}-\phi _{[k,\psi _{h}]^{0,1}}-\big ([\psi _{h},\phi _{k}]-\psi _{[h,\phi _{k}]^{0,1}}-\psi _{[k,\phi _{h}]^{0,1}}\big ).\]

Rappelons à présent que le choix d'une forme de contact \( \eta  \),
et donc d'un champ de Reeb \( R \), permet de définir une connexion
privilégiée sur \( X \), la connexion de Webster ; pour \( h\in H \),
la torsion \( T_{R,h}\in H \) et \( h\rightarrow T_{R,h} \) est
anti-\textbf{\( \mathbf{C} \)}-linéaire, donc définit un élément
\( T_{R,\cdot }\in \Omega ^{0,1}_{X}\otimes T^{1,0}_{X} \). On peut
aussi décomposer 

\[
T'X=\mathbf{C}R\oplus T^{1,0}_{X},\qquad \overline{\partial }=\left( \begin{array}{cc}
\overline{\partial } & \flat \\
T_{R,\cdot } & \overline{\partial }_{H}
\end{array}\right) ,\]
 où l'opérateur \( \flat :H\rightarrow H^{*} \) est défini par \[
\flat h=h\lrcorner d\eta ;\]
 on remarquera que \( \flat  \) envoie \( T^{1,0}_{X} \) sur \( \Omega ^{0,1}_{X} \),
et on peut l'étendre par produit extérieur en un opérateur \[
\flat :\Omega ^{0,i}_{X}\otimes T^{1,0}_{X}\rightarrow \Omega ^{0,i+1}_{X}.\]

Dans le cas particulier de \( X \), la structure de fibré en cercles
au-dessus de \( P^{1}\times P^{n-2} \) permet un choix priviliégié
de forme de contact \( \eta  \), à savoir la forme de connexion du
fibré \( \mathcal{O}(-1,1) \), de sorte que le vecteur de Reeb soit
l'action infinitésimale de \( S^{1} \) sur \( X \). Dans ce cas,
\( d\eta  \) provient de la base du fibré, plus précisément \( d(i\eta ) \)
est la courbure de \( \mathcal{O}(-1,1), \) de sorte qu'on a l'égalité
\[
-d\eta =F=-\omega _{P^{1}}+\omega _{P^{n-2}},\]
 où \( \omega _{P^{n}} \) est la forme de Fubini-Study de \( P^{n} \).
De plus, la torsion de Webster s'annule, et l'opérateur \( \overline{\partial } \)
du fibré \( T'X \) se réduit à \begin{equation}
\label{eq-dbar}
\overline{\partial }=\left( \begin{array}{cc}
\overline{\partial } & \flat \\
0 & \overline{\partial }_{H}
\end{array}\right) .
\end{equation}
 Par conséquent, la condition d'intégrabilité (\ref{con-int}) sur
\( \phi \in \Omega ^{0,1}_{X}\otimes T^{1,0}_{X} \) s'écrit \begin{equation}
\label{con-int2}
\overline{\partial }_{H}\phi +\frac{1}{2}[\phi ,\phi ]=0,\quad \flat \phi =0.
\end{equation}

\subsubsection*{Déformations}

Nous voulons montrer que l'espace des \( \phi  \) satisfaisant cette
condition d'intégrabilité est une variété de dimension infinie. Pour
cela, nous avons besoin d'un peu d'analyse. Soit \( \mathcal{H}^{\ell } \)
l'espace de Folland-Stein des fonctions sur \( X \) ayant \( \ell  \)
dérivées horizontales dans \( L^{2} \). Dans tout l'article, on choisira
\( \ell  \) suffisamment grand de sorte que \( \mathcal{H}^{\ell }\subset C^{0} \)
au moins, de sorte que \( \mathcal{H}^{\ell } \) soit une algèbre,
voir \cite[lemme 5.1]{BiqADE}.

Nous définissons alors l'espace \[
\mathcal{J}^{\ell }=\{\phi \in \mathcal{H}^{\ell }(\Omega ^{0,1}_{X}\otimes T^{1,0}_{X}),\, \overline{\partial }\phi +\frac{1}{2}[\phi ,\phi ]=0\}.\]
 Pour étudier \( \mathcal{J}^{\ell } \), nous devons faire une théorie
de déformation des structures CR de \( X \).

Remarquons que pour \( \phi \in \Omega ^{0,\cdot }_{X}\otimes T^{1,0}_{X} \),
on a \begin{equation}
\label{con-1}
\overline{\partial }(\phi \lrcorner d\eta )=(\overline{\partial }\phi )\lrcorner d\eta .
\end{equation}
 D'autre part, si \( \phi ,\psi \in \Omega ^{0,\cdot }_{X}\otimes T^{1,0}_{X} \)
sont dans le noyau de \( \flat  \), c'est-à-dire vérifient \[
\phi \lrcorner d\eta =\psi \lrcorner d\eta =0,\]
 alors, d'après (\ref{eq-dbar}), \( \overline{\partial }\phi  \)
et \( \overline{\partial }\psi  \) sont aussi éléments de \( \Omega ^{0,\cdot }_{X}\otimes T^{1,0}_{X} \),
de même que \( [\overline{\partial }\phi ,\psi ] \) et \( [\phi ,\overline{\partial }\psi ] \)
; à présent 

\[
\overline{\partial }[\phi ,\psi ]=[\overline{\partial }\phi ,\psi ]+(-1)^{\deg \phi }[\phi ,\overline{\partial }\psi ]\in \Omega ^{0,\cdot }_{X}\otimes T^{1,0}_{X},\]
 donc

\begin{equation}
\label{con-2}
[\phi ,\psi ]\lrcorner d\eta =0.
\end{equation}

Il résulte de (\ref{con-1}) et (\ref{con-2}) que \begin{equation}
\label{def-C}
\mathcal{C}^{\cdot }=(\ker \flat ,\overline{\partial }_{H})\subset \Omega ^{0,\cdot }_{X}\otimes T^{1,0}_{X}
\end{equation}
 est une algèbre différentielle graduée ; cette algèbre gouverne les
déformations de structures CR sur \( X \), nous allons donc pouvoir
lui appliquer la théorie classique de déformation pour étudier \( \mathcal{J}^{\ell } \).

Tout d'abord, nous avons besoin du calcul de certains groupes de cohomologie,
voir par exemple \cite[I.1]{OkoSchSpi80}.

\begin{lem}
\label{lem-HP}Pour \( 0<q<n \) et tout \( k \), on a \[
H^{q}(P^{n},\mathcal{O}(k))=H^{q}(P^{n},T_{P^{n}}(k))=0,\]
 à la seule exception de \( H^{n-1}(P^{n},T_{P^{n}}(-n-1))=\mathbf{C} \).

Enfin, on a \( H^{0}(P^{1},\mathcal{O}(k))=H^{1}(P^{1},\mathcal{O}(-k-2))=0 \)
si \( k<0 \).
\end{lem}
Remarquons que, puisque \( X \) est de signature \( (1,n-2) \),
l'opérateur \( \overline{\partial } \) sur \( X \) est hypoelliptique,
sauf en degrés \( 1 \) et \( n-2 \). Ainsi l'espace \( H^{2}(\mathcal{C}^{\cdot }) \)
est de dimension finie, sauf pour \( n=4 \) (en dimension \( 7 \)).
Plus précisément, on a le résultat suivant.

\begin{lem}
\label{lem-H2}En dimension \( 2n-1\neq 7 \), on a \( H^{2}(\mathcal{C}^{\cdot })=0 \).
\end{lem}
\begin{proof}
Les opérateurs \( \overline{\partial } \) et \( \flat  \) commutent
à l'action de \( S^{1} \), donc on peut se contenter de regarder
le problème en décomposant sous l'action de \( S^{1} \). Soit donc
\( \psi \in \Omega ^{0,2}_{X}\otimes T^{1,0}_{X} \) tel que \( \overline{\partial }_{H}\psi =0 \)
et \( \flat \psi =0 \), dans l'espace de poids \( k \) pour l'action
de \( S^{1} \). Cela signifie, en appelant \( t \) la section tautologique
du fibré \( \mathcal{O}(-1,1) \) sur \( X \), que \( \psi \otimes t^{-k} \)
descend en une section sur \( P^{1}\times P^{n-2} \) du fibré \( \Omega ^{0,2}\otimes T^{1,0}\otimes \mathcal{O}(k,-k) \).
Le problème se réduit donc, étant donnée sur \( P^{1}\times P^{n-2} \)
une section \( \psi  \) de \( \Omega ^{0,2}\otimes T^{1,0}\otimes \mathcal{O}(k,-k) \)
satisfaisant \begin{equation}
\label{con-psi}
\overline{\partial }\psi =0,\quad \psi \lrcorner F=0,
\end{equation}
 à trouver une section \( \phi  \) de \( \Omega ^{0,1}\otimes T^{1,0}\otimes \mathcal{O}(k,-k) \)
satisfaisant \begin{equation}
\label{con-phi}
\overline{\partial }\phi =\psi ,\quad \phi \lrcorner F=0.
\end{equation}

Or, pour tout entier \( k \), on a d'après le lemme \ref{lem-HP}
les égalités \begin{multline}
\label{for-H2a}
H^{2}(P^{1}\times P^{n-2},T_{P^{1}}^{1,0}\otimes \mathcal{O}(k,-k))\\=\bigoplus _{i=0,1}H^{i}(P^{1},\mathcal{O}(2+k))\otimes H^{2-i}(P^{n-2},\mathcal{O}(-k))=0,
\end{multline}
\begin{multline}
\label{for-H2b}
H^{2}(P^{1}\times P^{n-2},T_{P^{n-2}}^{1,0}\otimes \mathcal{O}(k,-k))\\=\bigoplus _{i=0,1}H^{i}(P^{1},\mathcal{O}(k))\otimes H^{2-i}(P^{n-2},T_{P^{n-2}}(-k))=0,
\end{multline}
 sauf dans le cas où \( n-2=3 \) où on obtient, puisque \( H^{2}(P^{3},T_{P^{3}}(-4))\neq 0 \),
\begin{equation}
\label{for-H2P1P3}
H^{2}(P^{1}\times P^{3},T_{P^{1}\times P^{3}}^{1,0}\otimes \mathcal{O}(4,-4))=H^{0}(P^{1},\mathcal{O}(4))H^{2}(P^{3},T_{P^{3}}(-4)).
\end{equation}

Dans l'immédiat, contentons-nous du cas \( n>5 \), donc il existe
\( \phi  \) telle que \( \overline{\partial }\phi =\psi  \). Reste
à modifier éventuellement \( \phi  \) de sorte de vérifier aussi
la seconde condition dans (\ref{con-phi}). Observons que, par (\ref{con-psi}),
\[
\overline{\partial }(\phi \lrcorner F)=(\overline{\partial }\phi )\lrcorner F=0.\]
 Comme \( H^{2}(P^{1}\times P^{n-2},\mathcal{O}(k,-k))=0 \) pour
\( n-2\geq 1 \), on en déduit l'existence d'une \( (0,1) \)-forme
\( \alpha  \), à valeurs dans \( \mathcal{O}(k,-k) \), telle que
\( \overline{\partial }\alpha =\phi \lrcorner F \). Il est alors
immédiat que \( \phi -\overline{\partial }\alpha  \) satisfait les
conditions (\ref{con-phi}).

Finalement, examinons plus précisément le cas \( n=5 \), où on semble
avoir l'espace de cohomologie non nul (\ref{for-H2P1P3}). Observons
tout d'abord que l'application \[
\lrcorner F:H^{2}(P^{3},T_{P^{3}}(-4))\rightarrow H^{3}(P^{3},T_{P^{3}}(-4))\]
 est manifestement un isomorphisme, entre deux espaces tous deux égaux
à \( \mathbf{C} \). Par conséquent, l'application \( \psi \rightarrow \psi \lrcorner F \)
définit un isomorphisme \begin{eqnarray*}
H^{2}(P^{1}\times P^{3},T^{1,0}_{P^{1}\times P^{3}}\otimes \mathcal{O}(4,-4)) & \rightarrow  & H^{3}(P^{1}\times P^{3},\mathcal{O}(4,-4))\\
 &  & =H^{0}(P^{1},\mathcal{O}(4))H^{3}(P^{3},T_{P^{3}}(-4)).
\end{eqnarray*}
 À présent, considérons une forme \( \psi  \) comme dans (\ref{con-psi}),
on déduit que la condition \( \psi \lrcorner F=0 \) impose l'annulation
de l'image de \( \psi  \) dans \( H^{2}(P^{1}\times P^{2},T^{1,0}_{P^{1}\times P^{3}}\otimes \mathcal{O}(4,-4)) \),
ce qui achève la démonstration du lemme.
\end{proof}
Le résultat essentiel de cette section est le suivant.

\begin{lem}
\label{lem-carte}L'espace \( \mathcal{J}^{\ell } \) est une sous-variété
hilbertienne de \( \mathcal{H}^{\ell }(\Omega ^{0,1}_{X}\otimes T^{1,0}_{X}) \),
d'espace tangent en \( 0 \) donné par les \( \dot{\phi }\in \mathcal{H}^{\ell }(\Omega ^{0,1}_{X}\otimes T^{1,0}_{X}) \)
tels que \( \overline{\partial }\dot{\phi }=0 \) (c'est-à-dire \( \overline{\partial }_{H}\dot{\phi }=0 \)
et \( \flat \dot{\phi }=0 \)). De plus, il existe une carte \[
\Phi :T_{0}\mathcal{J}^{\ell }\rightarrow \mathcal{J}^{\ell }\]
 telle que \( \Phi (\dot{\phi }) \) soit à coefficients de Fourier
positifs si et seulement si \( \dot{\phi } \) est à coefficients
de Fourier positifs.
\end{lem}
\begin{proof}
Nous faisons le cas de la dimension \( 2n-1\neq 7 \), le cas de la
dimension \( 7 \) est plus technique et fera l'objet de la section
\ref{sec-dim7}. 

À partir du lemme \ref{lem-H2}, la démonstration est classique. Choisissons
un inverse à droite \( P \) de l'opérateur \( \overline{\partial }_{H}:\mathcal{H}^{\ell }(\mathcal{C}^{1})\rightarrow \ker \overline{\partial }_{H}\subset \mathcal{H}^{\ell -1}(\mathcal{C}^{2}) \),
on peut choisir \( P \) commutant à l'action de \( S^{1} \). Étant
donné \( \dot{\phi } \) satisfaisant \( \overline{\partial }_{H}\dot{\phi }=0 \)
et \( \flat \dot{\phi }=0 \), nous construisons \[
\Phi (\dot{\phi })=\phi _{1}+\phi _{2}+\cdots ,\quad \phi _{1}=\dot{\phi },\quad \phi _{i}=-P\frac{1}{2}\sum _{1\leq j\leq i-1}[\phi _{j},\phi _{i-j}].\]
 Il est important d'écrire cette formule, car comme conséquence nous
voyons que si \( \dot{\phi } \) n'a que des coefficients de Fourier
positifs ou nuls, alors il en est de même pour \( \Phi (\dot{\phi }) \).

Réciproquement, à partir de \( \phi \in \mathcal{J}^{\ell } \), on
obtient \( \dot{\phi } \) par le choix \[
\dot{\phi }=\phi +P\frac{1}{2}[\phi ,\phi ].\]
 Là encore, \( \dot{\phi } \) n'a que des coefficients de Fourier
positifs s'il en était de même pour \( \phi  \).
\end{proof}

\section{\label{sec-2}Remplissage}

Ici nous rappelons certains résultats de \cite{BiqADE}, valables
en toute dimension, puis nous les appliquons à notre situation.

Habituellement, en théorie de la déformation, l'espace \( \mathcal{C}^{0} \)
est constitué des équivalences infinitésimales de la situation. Dans
la théorie décrite dans la section précédente, un élément de \( \mathcal{C}^{0} \)
est une section de \( \ker \flat \subset T^{1,0} \), donc en réalité
\( \mathcal{C}^{0}=0 \) et il n'y en a pas.

\subsubsection*{Action des contactomorphismes}

Bien entendu, il existe un groupe de jauge agissant sur les structures
CR, à savoir les contactomorphismes de \( X \). Rappelons que, d'après
\cite{Bla94} en dimension 3, et \cite[théorème 5.5]{BiqADE} en dimension
supérieure, il existe un groupe \( \mathcal{G}^{\ell +1} \) de contactomorphismes
de régularité \( \mathcal{H}^{\ell +1} \), paramétré près de l'origine
par les fonctions réelles de régularité \( \mathcal{H}^{\ell +2} \).
Infinitésimalement, une fonction réelle \( f \) correspond au contactomorphisme
infinitésimal \begin{equation}
\label{cont-inf}
fR-\sharp df,
\end{equation}
 où \( \sharp :H^{*}\rightarrow H \) est défini par \( (\sharp \alpha )\lrcorner d\eta =\alpha  \),
c'est-à-dire \( \sharp =\flat ^{-1} \). Le groupe \( \mathcal{G}^{\ell +1} \)
agit de manière lisse sur \( \mathcal{J}^{\ell } \) et l'action infinitésimale
à l'origine est donnée par \begin{equation}
\label{cont-act}
f\rightarrow \overline{\partial }_{H}\sharp \overline{\partial }f.
\end{equation}

Via la carte \( \Phi  \) définie par le lemme \ref{lem-carte}, on
peut voir \( \mathcal{G}^{\ell +1} \) comme agissant plutôt sur \( T_{0}\mathcal{J}^{\ell } \),
par \[
(\alpha ,\dot{\phi })\in \mathcal{G}^{\ell +1}\times T_{0}\mathcal{J}^{\ell }\longrightarrow \Phi ^{-1}\alpha \Phi (\dot{\phi }).\]
 L'action infinitésimale correspondante demeure donnée par l'opérateur
\( \overline{\partial }_{H}\sharp \overline{\partial } \), dont l'image
est fermée parce que l'opérateur \( \overline{\partial } \) est hypoelliptique
en degré 0 sur les variétés de signature \( (1,n-2) \). Cela implique
aussitôt le corollaire suivant \cite[section 6]{BiqADE}.

\begin{lem}
\label{lem-jauge}Supposons fixé un supplémentaire \( W \) de l'image
de \( \mathcal{H}^{\ell +2}(\mathbf{R}) \) par l'opérateur \( \overline{\partial }_{H}\sharp \overline{\partial } \)
dans \( T_{0}\mathcal{J}^{\ell } \), invariant sous l'action des
automorphismes CR de \( X \). Alors pour tout \( \dot{\phi }\in T_{0}\mathcal{J}^{\ell } \)
proche de \( 0 \), il existe un contactomorphisme \( \alpha  \),
proche de l'identité, tel que \( \Phi ^{-1}\alpha \Phi (\dot{\phi })\in W \)
; de plus \( \alpha  \) est unique modulo les automorphismes CR de
\( X \).\qed
\end{lem}

\subsubsection*{Condition de remplissage}

Rappelons que notre hypersurface réelle \( X=\partial N \) de \( P^{n} \)
est le fibré en cercles \( \mathcal{O}(-1,1) \) sur \( P^{1}\times P^{n-2} \).
Cependant, \( N \) diffère légèrement du fibré en disques \( \mathcal{O}(-1,1) \)
: d'après la formule (\ref{def-N}), celui-ci n'est autre que l'éclatement
\( \tilde{N} \) de \( N \) le long de \begin{equation}
\label{for-Pn2}
P^{n-2}=\{[0:0:z_{3}:\cdots :z_{n+1}]\}\subset N.
\end{equation}

Rappelons également que le fibré \( \Omega ^{0,1}_{X}\otimes T^{1,0}_{X} \)
provient du fibré \( \Omega ^{0,1}\otimes T^{1,0} \) sur la base
\( P^{1}\times P^{n-2} \) du fibré en cercles, ce qui permet de décomposer
le tenseur \( \phi \in \mathcal{J}^{\ell } \) définissant une structure
CR suivant ses coefficients de Fourier selon l'action de \( S^{1} \).

Dans cette situation, on a le résultat suivant \cite[théorème 4.1 et corollaire 4.2]{BiqADE}.

\begin{thm*}
Soit \( D\rightarrow Z \) un fibré holomorphe en disques, tel que
la forme de Levi du bord \( X=\partial D \) soit non dégénérée. Alors,
pour une structure CR intégrable sur \( X \), proche de la structure
standard, les deux propriétés suivantes sont équivalentes :
\begin{enumerate}
\item il existe un contactomorphisme de \( X \) qui envoie la structure
CR sur une structure dont les coefficients de Fourier strictement
négatifs s'annulent ;
\item la structure CR est le bord d'une déformation complexe de \( D \)
dont la section nulle \( Z \) demeure une sous-variété complexe.
\end{enumerate}
\end{thm*}
L'application dans notre situation donne la conséquence suivante.

\begin{cor}
Une structure CR sur \( X=\partial N \), proche de la structure standard,
est remplissable dans \( N \) si et seulement s'il existe un contactomorphisme
qui envoie la structure CR sur une structure dont les coefficients
de Fourier strictement négatifs s'annulent.
\end{cor}
\begin{proof}
Le théorème précédent dit que la condition sur les coefficients de
Fourier est la condition pour pouvoir remplir la structure CR dans
\( \tilde{N} \), de sorte que la section nulle \( P^{1}\times P^{n-2} \)
demeure complexe. Il reste à montrer que cela est équivalent à pouvoir
remplir la structure CR dans \( N \).

Si la structure est remplissable comme une déformation de \( N \),
alors le \( P^{n-2}\subset N \) donné par (\ref{for-Pn2}), de fibré
normal \( \mathcal{O}(1) \), persiste dans la déformation, et on
obtient une déformation de \( \tilde{N} \) en éclatant \( N \) le
long de \( P^{n-2} \).

Réciproquement, si une structure CR est le bord d'une déformation
de \( \tilde{N} \) pour laquelle la section nulle \( P^{1}\times P^{n-2} \)
reste complexe, alors la structure complexe de \( P^{1}\times P^{n-2} \)
reste la même puisque les déformations complexes de \( P^{1}\times P^{n-2} \)
sont triviales, et le fibré normal \( \mathcal{O}(-1,1) \) de \( P^{1}\times P^{n-2} \)
ne peut pas non plus changer ; cela signifie que l'on peut contracter
les \( P^{1} \) dans \( P^{1}\times P^{n-2} \) pour obtenir finalement
une déformation complexe de \( N \), dont la structure CR donnée
est le bord. 
\end{proof}
Autrement dit, pour savoir si une structure CR est remplissable, il
faut faire agir un contactomorphisme de sorte de tuer les coefficients
de Fourier négatifs. Compte tenu du lemme \ref{lem-jauge}, il suffit
de regarder cette question au niveau infinitésimal, sur l'image de
l'opérateur \( \overline{\partial }_{H}\sharp \overline{\partial } \).
Étant donnée une structure CR infinitésimale, on cherche donc une
fonction dont les coefficients négatifs tuent ceux de la structure
CR, les coefficients positifs étant fixés par la condition que la
fonction doit être réelle. Cette idée aboutit au corollaire suivant
\cite[lemme 6.4]{BiqADE}. Notons \( \pi _{-} \) la projection sur
les coefficients de Fourier négatifs ou nuls.

\begin{lem}
\label{lem-coulomb}Choisissons comme supplémentaire de l'image de
\( \mathcal{H}^{\ell +2}(\mathbf{R}) \) par \( \overline{\partial }_{H}\sharp \overline{\partial } \)
l'espace \[
W=(\im \overline{\partial }_{H}\sharp \overline{\partial }\pi _{-})^{\perp },\]
 alors une structure CR intégrable \( \phi =\Phi (\dot{\phi })\in \mathcal{J}^{\ell } \),
proche de 0, est remplissable si et seulement si dans la jauge \( \Phi ^{-1}\alpha \Phi (\dot{\phi })\in W \)
construite par le lemme \ref{lem-jauge}, les coefficients de Fourier
strictement négatifs s'annulent.\qed
\end{lem}

\subsubsection*{Calcul des obstructions}

Le lemme \ref{lem-coulomb} n'est utile que si l'on sait déterminer
l'image de l'opérateur \( \overline{\partial }_{H}\sharp \overline{\partial } \).
Il suffit de regarder ce qui se passe dans chaque espace de poids
\( k \) sous l'action de \( S^{1} \).

Commençons par remarquer que l'on peut définir un complexe étendu
\( \tilde{\mathcal{C}} \) en remplaçant l'espace trivial \( \mathcal{C}^{0} \)
du complexe (\ref{def-C}) par un espace non trivial \begin{equation}
\label{ref-Ctilde}
\tilde{\mathcal{C}}^{0}=\{\upsilon \in \mathcal{H}^{\ell +1}(T^{1,0}_{X}),\, (\overline{\partial }_{H}\upsilon )\lrcorner F=0\}.
\end{equation}
 La condition sur \( \upsilon  \) est équivalente à \( \overline{\partial }(\upsilon \lrcorner F)=0 \),
et des exemples de tels \( \upsilon  \) sont fournis comme \( \sharp \overline{\partial }f \),
où \( f \) est une fonction complexe ; comme \( \upsilon \lrcorner F \)
est une (0,1)-forme, il y a beaucoup d'autres \( \upsilon  \), puisque
l'opérateur \( \overline{\partial } \) n'est pas hypoelliptique en
ce degré.

\begin{lem}
\label{lem-H1}En poids \( k\leq 0 \), le \( H^{1} \) du complexe
étendu \( \tilde{\mathcal{C}} \) s'identifie à \begin{multline*}
H^{1}(P^{1}\times P^{n-2},T^{1,0}_{P^{1}\times P^{n-2}}\otimes \mathcal{O}(k,-k)) =\\ H^{1}(P^{1},T^{1,0}_{P^{1}}(k))H^{0}(P^{n-2},\mathcal{O}(-k))
  \oplus H^{1}(P^{1},\mathcal{O}(k))H^{0}(P^{n-2},T^{1,0}_{P^{n-2}}(-k)).
\end{multline*}

\end{lem}
\begin{proof}
On a clairement une application injective \[
H^{1}(\tilde{\mathcal{C}})\rightarrow H^{1}(P^{1}\times P^{n-2},T^{1,0}_{P^{1}\times P^{n-2}}\otimes \mathcal{O}(k,-k)).\]
 Réciproquement, il faut vérifier que tout élément du second espace
de cohomologie provient bien de \( \tilde{\mathcal{C}} \) : nous
le faisons pour les deux morceaux écrits dans l'énoncé du lemme.

Le premier morceau est le plus simple : un élément dans \[
H^{1}(P^{1},T^{1,0}_{P^{1}}(k))H^{0}(P^{n-2},\mathcal{O}(-k))\]
 peut être représenté par une section de \( \Omega ^{0,1}_{P^{1}}\otimes T^{1,0}_{P^{1}}\otimes \mathcal{O}(k,-k) \),
qui est dans le noyau de \( \flat  \).

Le second morceau est plus compliqué : un élément dans \[
H^{1}(P^{1},\mathcal{O}(k))H^{0}(P^{n-2},T^{1,0}_{P^{n-2}}(-k))\]
 est représenté par une section \( \phi \in \Omega ^{0,1}_{P^{1}}\otimes T^{1,0}_{P^{n-2}}\otimes \mathcal{O}(k,-k) \),
holomorphe le long de \( P^{n-2} \) ; on remplace \( \phi  \) par
\( \phi +\overline{\partial }\upsilon  \), où \( \upsilon \in T^{1,0}_{P^{1}}\otimes \mathcal{O}(k,-k) \)
satisfait l'équation \[
\overline{\partial }\upsilon \lrcorner F=-\phi \lrcorner F.\]
 Cette équation s'écrit encore \( \overline{\partial }(\upsilon \lrcorner F)=-\phi \lrcorner F \),
où \( \upsilon \lrcorner F\in \Omega ^{0,1}_{P^{1}}\otimes \mathcal{O}(k,-k) \)
et \( \phi \lrcorner F\in \Omega ^{0,1}_{P^{1}}\otimes \Omega ^{0,1}_{P^{n-2}}\otimes \mathcal{O}(k,-k) \),
et l'existence de la solution provient de l'annulation \( H^{1}(P^{n-2},\mathcal{O}(-k))=0 \).
\end{proof}
À présent, nous pouvons calculer le supplémentaire \( W \) du lemme
\ref{lem-coulomb}. Fixons un poids \( k\leq 0 \) de l'action de
\( S^{1} \), soit \( \phi \in T_{0}\mathcal{J}^{\ell } \) en poids
\( k \), nous cherchons une section \( f \) de \( \mathcal{O}(k,-k) \)
telle que \( \phi =\overline{\partial }_{H}\sharp \overline{\partial }f \).
La première obstruction est évidemment \[
[\phi ]=0\in H^{1}(\tilde{\mathcal{C}});\]
 si cela est le cas, il existe \( \upsilon \in \tilde{\mathcal{C}}^{0} \)
tel que \( \overline{\partial }_{H}\upsilon =\phi  \). Observons
que \begin{equation}
\label{for-H0}
H^{0}(P^{1}\times P^{n-2},T^{1,0}\otimes \mathcal{O}(k,-k))=H^{0}(P^{1},T^{1,0}_{P^{1}}(k))H^{0}(P^{n-2},\mathcal{O}(-k))
\end{equation}
 est trivial, sauf en degrés \( k=0 \), \( -1 \) ou \( -2 \). Par
conséquent, \( \upsilon  \) est unique, sauf en degrés \( 0 \),
\( -1 \) ou \( -2 \) où l'ambiguïté est mesurée par (\ref{for-H0}). 

Maintenant, \( \flat \upsilon  \) est une (0,1)-forme satisfaisant
d'après (\ref{con-1}) \[
\overline{\partial }\flat \upsilon =\flat \overline{\partial }\upsilon =\flat \phi =0,\]
 et nous cherchons \( f \) telle que \( \overline{\partial }f=\upsilon  \).
Observons que l'espace de cohomologie correspondant, \[
H^{1}(P^{1}\times P^{n-2},\mathcal{O}(k,-k))=H^{1}(P^{1},\mathcal{O}(k))H^{0}(P^{n-2},\mathcal{O}(-k)),\]
 est non trivial en général, donc il y a une seconde obstruction.
Pour \( k=0 \) ou \( -1 \), cette obstruction s'annule. Pour \( k=-2 \),
l'obstruction correspond exactement à la liberté (\ref{for-H0}) pour
\( \upsilon  \), donc il n'y a pas d'obstruction supplémentaire.
En revanche, pour \( k<-2 \), compte tenu de l'unicité de \( \upsilon  \),
l'obstruction correspondante pour \( \phi  \) est l'espace \( (\overline{\partial }_{P^{1}}\sharp H^{1}(P^{1},\mathcal{O}(k)))H^{0}(P^{n-2},\mathcal{O}(-k)) \),
où cette notation signifie qu'un sous-espace représentant \( H^{1}(P^{1},\mathcal{O}(k)) \)
a été choisi. Nous résumons la discussion précédente dans le résultat
suivant.

\begin{lem}
\label{lem-obstructions}Pour \( k\leq 0 \), l'orthogonal de l'image
de \( \overline{\partial }_{H}\sharp \overline{\partial } \) en poids
\( k \) s'identifie à \begin{eqnarray*}
W_{k} & = & H^{1}(P^{1},\mathcal{O}(k))H^{0}(P^{n-2},T^{1,0}_{P^{n-2}}(-k))\\
 &  & \oplus \left( H^{1}(P^{1},T^{1,0}_{P^{1}}(k))\oplus \overline{\partial }_{P^{1}}\sharp H^{1}(P^{1},\mathcal{O}(k))\right) H^{0}(P^{n-2},\mathcal{O}(-k)).
\end{eqnarray*}
\qed
\end{lem}
Dans ce lemme il faut noter que pour \( k=-1 \) et \( -2 \), on
a $$ \overline{\partial }_{P^{1}}\sharp H^{1}(P^{1},\mathcal{O}(k))=0   .$$
Puisque \( W_{k}\neq 0 \) pour tout \( k<-1 \), on déduit finalement
de ce calcul le théorème \ref{th-def}.

\begin{rem}
\label{rem-kir}Par le théorème de Kiremidjian (voir l'introduction),
une petite déformation de \( X \) est toujours remplissable de l'autre
côté, c'est-à-dire du côté de \( M \) défini par (\ref{def-M}).
La méthode utilisée dans cet article permet de donner une autre démonstration
de ce fait, plus simple que celle de Kiremidjian qui utilise le théorème
de Nash-Moser. En effet, un calcul similaire à celui de la démonstration
précédente montre que si \( n>3 \), alors \( W_{k}=0 \) pour \( k\geq 0 \),
si bien qu'une structure CR sur \( X \), proche de la structure standard,
peut toujours être ramenée par un contactomorphisme à une structure
dont les coefficients de Fourier positifs s'annulent, donc remplissable
du côté de \( M \).
\end{rem}

\begin{rem}
\label{rem-0}Le cas des coefficients invariants sous \( S^{1} \)
(\( k=0 \)) est un peu particulier, puisque la condition que la fonction
\( f \) qui paramètre un contactomorphisme soit réelle donne une
restriction sur \( f \) en poids \( 0 \). Comme on va le voir dans
la section \ref{sec-3}, les fonctions \( f \) complexes représentent
les déformations de \( X \) dans \( P^{n} \), et le fait que \( W_{0}=0 \)
indique que les structures CR \( S^{1} \)-invariantes sur \( X \)
sont obtenues par déformation (également \( S^{1} \)-invariante)
de \( X \) dans \( P^{n} \).
\end{rem}

\section{\label{sec-3}Stabilité}

\subsubsection*{Aspects formels}

Ici nous ne donnons pas tous les détails, car nous démontrerons le
résultat \ref{th-stab} par une autre méthode.

On peut obtenir d'autres structures CR sur \( X \), en déformant
l'hypersurface réelle \( X \) dans \( P^{n} \). Si on a une application
\( F:X\rightarrow P^{n} \), on en déduit sur \( X \) une nouvelle
structure CR. Si on s'intéresse seulement aux petites déformations
de la sous-variété initiale \( X \), alors la structure de contact
sous-jacente à la structure CR demeure équivalente à la structure
de contact initiale, via un difféomorphisme \( \alpha  \) de \( X \).
Par conséquent, quitte à remplacer \( F \) par \( F\circ \alpha  \),
on peut se limiter à regarder les applications \( F:X\rightarrow P^{n} \)
telles que la structure CR induite sur \( X \) conserve la même structure
de contact sous-jacente. Cette condition sur \( F \) s'écrit \begin{equation}
\label{con-F}
Id_{x}F(H_{x})\subset d_{x}F(H_{x}),
\end{equation}
où \( I \) est la structure complexe de \( V \).

Introduisons l'élément \( \eta ^{c}=\eta +iI\eta  \) de \( \Omega ^{1,0}_{V}|_{X}=(T'X)^{*} \).
Un calcul sans difficulté montre qu'infinitésimalement, la condition
(\ref{con-F}) devient, pour \( \dot{F} \) section de \( TP^{n}|_{X} \)
et \( h\in H \), \[
\eta ^{c}(I\overline{\partial }_{h}\dot{F})=0.\]
 Écrivant \( \dot{F}^{1,0}=fR+h \), cette condition se récrit \[
\overline{\partial }f+h\lrcorner d\eta =0.\]
 On observe donc qu'infinitésimalement, les applications \( F \)
satisfaisant la condition (\ref{con-F}) sont données par les \begin{equation}
\label{defo-inf}
f-\sharp \overline{\partial }f,
\end{equation}
 où \( f \) est une fonction \emph{complexe} sur \( X \).

Comparant (\ref{defo-inf}) à (\ref{cont-inf}), il est clair que
l'espace des applications infinitésimales satisfaisant (\ref{con-F})
est une complexification de l'espace des contactomorphismes infinitésimaux
; de plus il est facile de vérifier que l'action infinitésimale de
\( f-\sharp \overline{\partial }f \) sur la structure CR est \[
f\rightarrow \overline{\partial }_{H}\sharp \overline{\partial }f,\]
 soit la complexification de l'action (\ref{cont-act}). 

Pour résumer, l'action des déformations de \( X \) dans \( P^{n} \),
préservant la structure de contact de \( X \), sur les structures
CR de \( X \), apparaît comme complexification de l'action des contactomorphismes.
Cette remarque est utile, mais elle demeure à un niveau formel : en
effet, pour la rendre rigoureuse, il faut construire une paramétrisation
des applications de régularité \( \mathcal{H}^{\ell +1} \) satisfaisant
(\ref{con-F}) par des fonctions de régularité \( \mathcal{H}^{\ell +2} \),
ce qui semble difficile quand l'opérateur \( \overline{\partial } \)
n'est pas hypoelliptique en degré 1. Comme conséquence d'une telle
construction, on aurait une démonstration plus simple du théorème
de Hamilton.

Poursuivant ce point de vue formel, on voit que, de manière analogue
au lemme \ref{lem-coulomb}, une structure CR \( \Phi (\dot{\phi }) \)
sur \( X \) provient d'une déformation de \( X \) dans \( P^{n} \)
si \( \dot{\phi } \) est dans l'image de \( \mathcal{H}^{\ell +2}(\mathbf{C}) \)
par l'opérateur \( \overline{\partial }_{H}\sharp \overline{\partial } \).
En poids \( k<0 \), cela signifie que la structure CR est remplissable
du côté de \( N \) ; en poids \( k>0 \), cela signifie que la structure
CR est remplissable de l'autre côté, voir la remarque \ref{rem-kir}
; pour le poids \( k=0 \), voir la remarque \ref{rem-0}.

Une structure CR sur \( X \), proche de la structure standard, et
stable (au sens donné dans l'introduction), est évidemment remplissable
des deux côtés. La réciproque est vraie, car si on a un remplissage
des deux côtés, alors en recollant le long de \( X \), on obtient
une variété complexe compacte qui est une déformation de \( P^{n} \),
donc exactement \( P^{n} \). Cette observation aboutit à la démonstration
suivante.

\subsubsection*{Démonstration du corollaire \ref{th-stab}}

Supposons donnée sur \( X^{2n-1} \) (\( n>3 \)) une structure CR
intégrable, proche de la structure standard. Par le théorème de Kiremidjian
(ou la remarque \ref{rem-kir}), \( X \) est toujours remplissable
du côté de \( M \). 

Supposons ce remplissage stable, c'est-à-dire induit par une déformation
de \( X \) dans \( P^{n} \). Alors il est clair que \( X \) est
aussi remplissable du côté de \( N \).

Réciproquement, supposons \( X \) remplissable du côté de \( N \).
Par le théorème de Hamilton, ce remplissage est induit par une déformation
de \( X \) dans \( P^{n} \), donc il y a stabilité pour le remplissage
de \( X \) du côté de \( M \) aussi.

Il en résulte que le remplissage de \( X \) du côté de \( M \) est
stable si et seulement si \( X \) est remplissable du côté de \( N \),
c'est-à-dire vérifie les conditions du théorème \ref{th-def}. En
particulier les obstructions sont les mêmes que pour ce problème,
donc données par le lemme \ref{lem-obstructions}.\qed 

\begin{rem}
Dans \cite{CatLem92}, Catlin et Lempert construisent des exemples
\( X^{2k-1} \) de variétés CR, strictement pseudoconvexes, plongées
dans \( P^{n} \) (\( n>k \)), dont les déformations CR ne sont pas
induites par une déformation du plongement. Il est possible d'analyser
ce problème en termes similaires aux déformations qui apparaissent
dans cette section---l'absence de déformation du plongement provient
alors d'une obstruction cohomologique explicite.
\end{rem}

\section{Le cas de la dimension 7\label{sec-dim7}}

En dimension 7, la démonstration du lemme \ref{lem-carte} est compliquée
par le fait que le \( H^{2} \) du complexe \( \mathcal{C} \) est
non nul, et même de dimension infinie car \( \overline{\partial } \)
n'est pas hypoelliptique en degré 2. On explique ici comment remédier
à ce problème.

Fixons donc \( 2n-1=7 \). Comme pour \( k\geq 0 \), \( k\neq 3 \),
on a par le lemme \ref{lem-HP} \[
H^{1}(P^{1}\times P^{2},T^{1,0}\otimes \mathcal{O}(k,-k))=0,\]
 on déduit, comme dans le lemme \ref{lem-H1}, qu'en un poids \( k\geq 0 \),
\( k\neq 3 \), la cohomologie \( H^{1}(\tilde{\mathcal{C}}) \) s'annule.
Pour \( k=3 \), on a \[
H^{1}(P^{1}\times P^{2},T^{1,0}\otimes \mathcal{O}(3,-3))=H^{0}(P^{1},\mathcal{O}(3))H^{1}(P^{2},T^{1,0}P^{2}(-3)),\]
 mais par le même raisonnement tenu pour traiter (\ref{for-H2P1P3}),
on montre que malgré tout la cohomologie \( H^{1}(\tilde{\mathcal{C}}) \)
s'annule encore. Il en résulte finalement que \( H^{1}(\tilde{\mathcal{C}}) \)
est concentré sur les poids \( k<0 \).

D'un autre côté, en reprenant les formules (\ref{for-H2a}) et (\ref{for-H2b}),
on voit que pour \( k\leq 0 \) \[
H^{2}(P^{1}\times P^{2},T^{1,0}\otimes \mathcal{O}(k,-k))=0,\]
 donc au contraire \( H^{2}(\tilde{\mathcal{C}}) \) est concentré
sur les poids \( k>0 \).

En particulier, on voit que l'image du crochet \( H^{1}(\tilde{\mathcal{C}})H^{1}(\tilde{\mathcal{C}})\rightarrow H^{2}(\tilde{\mathcal{C}}) \)
est concentré en degrés négatifs, donc ne peut être que nulle dans
\( H^{2}(\tilde{\mathcal{C}}) \) :

\begin{lem}
L'application \[
H^{1}(\tilde{\mathcal{C}})H^{1}(\tilde{\mathcal{C}})\rightarrow H^{2}(\tilde{\mathcal{C}})\]
 induite par le crochet est identiquement nulle.\qed
\end{lem}
À partir de là, il est assez facile de déduire le lemme \ref{lem-carte}.
En effet, tout élément de \( T_{0}\mathcal{J}^{\ell } \) peut être
représenté par \[
\dot{\phi }+\overline{\partial }\upsilon ,\]
 où \( \dot{\phi } \) est un représentant à coefficients de Fourier
négatifs de la classe induite dans \( H^{1}(\tilde{\mathcal{C}}) \),
et \( \upsilon \in \tilde{\mathcal{C}}^{0} \). Pour \( \dot{\phi } \)
seul, le procédé de génération de \( \Phi (\dot{\phi }) \), comme
dans la démonstration du lemme \ref{lem-carte}, par \[
\Phi (\dot{\phi })=\phi _{1}+\phi _{2}+\cdots ,\quad \phi _{1}=\dot{\phi },\quad \phi _{i}=-P\frac{1}{2}\sum _{1\leq j\leq i-1}[\phi _{j},\phi _{i-j}],\]
 fonctionne, puisqu'on reste avec des coefficients de Fourier négatifs,
alors que \( H^{2} \) est concentré en poids positifs. Par des arguments
généraux, cela signifie que la résolution pour \( \dot{\phi }+\overline{\partial }\upsilon  \)
est également possible.

\providecommand{\bysame}{\leavevmode ---\ }
\providecommand{\og}{``}
\providecommand{\fg}{''}
\providecommand{\smfandname}{et}
\providecommand{\smfedsname}{\'eds.}
\providecommand{\smfedname}{\'ed.}
\providecommand{\smfmastersthesisname}{M\'emoire}
\providecommand{\smfphdthesisname}{Th\`ese}

\smallskip{}
{\raggedleft \noun{Irma}, Université Louis Pasteur et \noun{Cnrs},\par}

{\raggedleft 7 rue René Descartes, F-67084 Strasbourg Cedex\par}

\smallskip{}
\raggedleft \texttt{olivier.biquard@math.u-strasbg.fr}

\begin{thebibliography}{Ham77}

\bibitem[BE90]{BurEps90}
{\scshape D.~M. Burns {\normalfont \smfandname} C.~L. Epstein} -- {\og
  Embeddability for three-dimensional {C}{R}-manifolds\fg}, \emph{J. Amer.
  Math. Soc.} \textbf{3} (1990), no.~4, p.~809--841.

\bibitem[Biq02]{BiqADE}
{\scshape O.~Biquard} -- {\og M{\'e}triques autoduales sur la boule\fg},
  \emph{Invent. Math.} (2002).

\bibitem[Bla94]{Bla94}
{\scshape J.~S. Bland} -- {\og Contact geometry and {C}{R} structures on
  ${S}\sp 3$\fg}, \emph{Acta Math.} \textbf{172} (1994), no.~1, p.~1--49.

\bibitem[CL92]{CatLem92}
{\scshape D.~Catlin {\normalfont \smfandname} L.~Lempert} -- {\og {A note on
  the instability of embeddings of Cauchy-Riemann manifolds.}\fg}, \emph{J.
  Geom. Anal.} \textbf{2} (1992), no.~2, p.~99--104.

\bibitem[Eps92]{Eps92}
{\scshape C.~L. Epstein} -- {\og C{R}-structures on three-dimensional circle
  bundles\fg}, \emph{Invent. Math.} \textbf{109} (1992), no.~2, p.~351--403.

\bibitem[Ham77]{Ham77}
{\scshape R.~S. Hamilton} -- {\og Deformation of complex structures on
  manifolds with boundary. {I}. {T}he stable case\fg}, \emph{J. Differential
  Geom.} \textbf{12} (1977), no.~1, p.~1--45.

\bibitem[Kir79]{Kir79}
{\scshape G.~K. Kiremidjian} -- {\og A direct extension method for {C}{R}
  structures\fg}, \emph{Math. Ann.} \textbf{242} (1979), no.~1, p.~1--19.

\bibitem[Lem92]{Lem92}
{\scshape L.~Lempert} -- {\og On three-dimensional {C}auchy-{R}iemann
  manifolds\fg}, \emph{J. Amer. Math. Soc.} \textbf{5} (1992), no.~4,
  p.~923--969.

\bibitem[OSS80]{OkoSchSpi80}
{\scshape C.~Okonek, M.~Schneider {\normalfont \smfandname} H.~Spindler} --
  \emph{Vector bundles on complex projective spaces}, Birkh\"auser Boston,
  Mass., 1980.

\end{thebibliography}
\end{document}